\documentclass[11pt, a4paper]{article}

\usepackage{lipsum}
\usepackage{authblk}
\usepackage{bm}
\usepackage{amssymb, amsmath}
\usepackage{color}
\usepackage[amsmath,thmmarks]{ntheorem}
\numberwithin{equation}{section}
\newcommand{\fnum}{\mathbb{F}}
\DeclareMathOperator{\tr}{tr}
\DeclareMathOperator{\norm}{N}
\DeclareMathOperator{\moda}{mod}
\newcommand{\bases}[1]{\langle #1 \rangle}
\newtheorem{theorem}{Theorem}[section]

\newtheorem{lemma}[theorem]{Lemma}
\newtheorem{corollary}[theorem]{Corollary}

\theorembodyfont{\rmfamily}

\newtheorem{definition}{Definition}[section]

\theoremstyle{nonumberplain}
\theoremsymbol{$\square$}
\newtheorem{proof}{Proof}

\begin{document}

\title{Davenport--Hasse's Theorem for Polynomial Gauss Sums over Finite Fields}

\date{}

\author{Zhiyong Zheng\footnote{This work was partially supported by the ``973'' project 2013CB834205.}}
\affil{School of Mathematics and Systems Science, \authorcr
Beihang University, Beijing, P.R. China \authorcr
zhengzhiyong@buaa.edu.cn \vspace{4mm}}

\maketitle
\begin{abstract}
In this paper, we study the polynomial Gauss sums over finite fields, and present an analogue of Davenport--Hasse's theorem for the entire polynomial Gauss sums, which is a generalization of the previous result obtained by Hayes. 
\end{abstract}

\noindent{\small {\bf 2010 Mathematics Subject Classification: }Primary 11T55, 11T24, 11L05}

\noindent{\small {\bf Key words: }Polynomial Gauss sums, Finite fields, Davenport--Hasse's Theorem}

\section{Introduction}
Let $\fnum_q$ be a finite field with $q=p^l$ elements, where $p$ is a prime number. Let $\fnum_{q^n}$ be a finite extension of $\fnum_q$ of degree $n$, and $\sigma$ be the Frobinus on $\fnum_{q^n}$, given by $\sigma(a) = a^q$ for any element $a$ in $\fnum_{q^n}$. We have $\sigma^n =1$, and $\sigma$ generates the Galois group of $\fnum_{q^n}/\fnum_q$. The relative trace $\tr(a)$ and the norm $\norm(a)$ of an element $a$ in $\fnum_{q^n}$ are defined by 
\begin{equation}
  \label{eq:traceNorm}
\tr(a) = \sum_{i=1}^n \sigma^i(a), \quad \norm(a) = \prod_{i=1}^n \sigma^i(a)  
\end{equation}
respectively. Let $\psi$ be a (complex-valued) character of the additive group of $\fnum_q$, and $\chi$ a character of the multiplicative group $\fnum_q^*$ of $\fnum_q$. The Gauss sums of $\fnum_q$ are defined by 
\begin{equation}
  \label{eq:gaussSum}
  \tau(\chi, \psi) = \sum_{a\in \fnum_q^*}\chi(a)\psi(a).
\end{equation}
If we set for every $a$ in $\fnum_{q^n}$ that 
\begin{equation}
  \label{eq:character}
\psi^{(n)}(a) = \psi(\tr(a)), \quad \chi^{(n)}(a) = \chi(\norm(a)),  
\end{equation}
then the function $\psi^{(n)}$ will be a character of additive group, and function $\chi^{(n)}$ a character of the multiplicative group of $\fnum_{q^n}$. In particular, if $a\in \fnum_q$, then we have
\begin{equation}
  \label{eq:identity}
  \psi^{(n)}(a) = \psi^n(a), \quad \chi^{(n)}(a) =\chi^n(a).
\end{equation}
These characters $\psi^{(n)}$ and $\chi^{(n)}$ define a generalized Gauss sums $\tau(\chi^{(n)}, \psi^{(n)})$ on $\fnum_{q^n}$. Davenport and Hasse in \cite{DH1934D} proved the following remarkable theorem (also see \cite{IR1982a} and \cite{S1937r}) that

\begin{theorem}[Davenport--Hasse]\label{thm:dh} If both of $\chi$ and $\psi$ are not principal, then
  \begin{equation}
    \label{eq:theorem}
-\tau(\chi^{(n)},\psi^{(n)})  = (-\tau(\chi, \psi))^n.    
  \end{equation}
\end{theorem}

To generalize this famous theorem to the polynomial Gauss sums, let $\fnum_q[x]$ and $\fnum_{q^n}[x]$ be the polynomial rings over $\fnum_q$ and $\fnum_{q^n}$ respectively. The Frobinus $\sigma$ of $\fnum_{q^n}$ may be extended to $\fnum_{q^n}[x]$ in the following way: If $A = a_mx^m + \cdots + a_1x+a_0 \in \fnum_{q^n}[x]$, then we define
\begin{equation}
  \label{eq:sigmaA}
  \sigma(A) = \sigma(a_m)x^m + \cdots + \sigma(a_1)x + \sigma(a_0),
\end{equation}
which clearly is a $\fnum_q[x]$-automorphism of $\fnum_{q^n}[x]$. The relative trace and norm may be extended to $\fnum_{q^n}[x]$ by 
\begin{equation}
  \label{eq:traceNormGen}
  \tr(A) = \sum_{i=1}^n \sigma^i(A), \quad \norm(A) = \prod_{i=1}^n \sigma^i(A).
\end{equation}
Then $\tr(A)$ is an additive and $\norm(A)$ a multiplicative function from $\fnum_{q^n}[x]$ to $\fnum_q[x]$. 

Let $H$ be a fixed but arbitrary polynomial in $\fnum_q[x]$ with degree $m$, $\psi$ be a (complex-valued) character of additive group of the residue class ring $\fnum_q[x]/\bases{H}$. We may understand that $\psi$ is a complex-valued function defined on $\fnum_q[x]$ such that
\begin{equation}
  \label{eq:psi}
  \psi(A+B) = \psi(A)\cdot \psi(B), \mbox{and }\psi(A)=\psi(B), \mbox{ if } A \equiv B (\moda H)
\end{equation}
for any two polynomials $A$ and $B$ in $\fnum_q[x]$. According to Hayes \cite{H1966a}, we call $\psi$ an additive character modulo $H$ on $\fnum_q[x]$. For example, $\psi_0$ is the principal additive character modulo $H$, where $\psi_0(A)=1$ for all the polynomials $A$ in $\fnum_q[x]$. 

Let $\chi$ be a (complex-valued) character of the multiplicative group of the reduced residue of $\fnum_q[x]/\bases{H}$, $\chi$ may be also understood as a complex-valued function of $\fnum_q[x]$, such that $\chi(A) = 0$ if $(A,H)>1$ and
\begin{equation}
  \label{eq:chi}
  \chi(AB) = \chi(A)\cdot \chi(B), \mbox{and }\chi(A)=\chi(B), \mbox{ if } A \equiv B (\moda H).
\end{equation}
We also call $\chi$ a multiplicative character modulo $H$ on $\fnum_q[x]$. Especially, $\chi_0$ is the principal multiplicative character modulo $H$, where $\chi_0(A) = 1$ for all of $A$ in $\fnum_q[x]$ with $(A,H)=1$.

With the above notations, we define a polynomial Gauss sum $G(\chi, \psi)$ modulo $H$ on $\fnum_q[x]$ as follows
\begin{equation}
  \label{eq:gaussSum}
  G(\chi, \psi) = \sum_{D\moda H}\chi(D)\psi(D),
\end{equation}
where the summation extending over a complete residue system of modulo $H$ in $\fnum_q[x]$. 

For a polynomial $H$ in $\fnum_q[x]$ and, therefore, also a polynomial in $\fnum_{q^n}[x]$, to define a Gauss sum modulo $H$ on $\fnum_{q^n}[x]$, for any $A\in \fnum_{q^n}[x]$, we set $\psi^{(n)}(A)$ and $\chi^{(n)}(A)$ by 
$$\psi^{(n)}(A) = \psi(\tr(A)), \quad \chi^{(n)}(A) = \chi(\norm(A)).$$
It is easy to verify that $\psi^{(n)}$ is an additive and $\chi^{(n)}$ a multiplicative character modulo $H$ on $\fnum_{q^n}[x]$, thus we may use these characters to define a polynomial Gauss sum $G(\chi^{(n)}, \psi^{(n)})$ modulo $H$ on $\fnum_{q^n}[x]$. 

The most interesting question is that, is there an analogue of Davenport--Hasse's theorem for the polynomial Gauss sums? Hayes \cite[Theorem~2.2]{H1966a} shows that an analogue of Davenport--Hasse's theorem for a special additive character $\psi = E$, which character was essentially introduced by Carlitz (see \cite{C1947t}). To state Hayes' result, for any polynomial $A$ in $\fnum_{q}[x]$, let 
$$A\equiv a_{m-1}x^{m-1} + \cdots + a_1x + a_0 (\moda H),$$
where $m=\deg(H)$. We set an additive function modulo $H$ on $\fnum_q[x]$ by $t(A) = a_{m-1}$. By the definition, we have immediately that for any $A, B\in \fnum_q[x]$, $t(A+B) = t(A) + t(B)$, and $t(A)= t(B)$ whenever $A\equiv B (\moda H)$, in particular, $t(A)=0$ if $H|A$. To generalize that $t$-function modulo $H$ on $\fnum_q[x]$, for a given polynomial $G$ in $\fnum_q[x]$, let
\begin{equation}
  \label{eq:tG}
  t_G(A) = t(GA).
\end{equation}
Clearly, $t_G$ is an additive function modulo $H$ on $F_q[x]$, that is 
\begin{equation}
  \label{eq:tGAdd}
  t_G(A+B) = t_G(A)+t_G(B), \mbox{ and }t_G(A) = t_G(B), \mbox{ if } A\equiv B (\moda H).
\end{equation}

Now, let $\lambda$ be a fixed non-principal additive character on $\fnum_q$, for example, $\lambda(a) = e(\frac{2\pi i \tr(a)}{p})$, for $a$ in $\fnum_q$, we define the function $E_{\lambda}(G, H)$ on $\fnum_q[x]$ by 
\begin{equation}
  \label{eq:Elambda}
  E_{\lambda}(G, H)(A) = \lambda(t_G(A)).
\end{equation}
It is easily seen that $E_{\lambda}(G, H)$ is an additive character modulo $H$ on $\fnum_q[x]$. If we set $G=1$, a constant polynomial, and let $E=E_{\lambda} = E_{\lambda}(1, H)$, then Hayes \cite{H1966a} shows that an analogue of Davenport--Hasse's theorem for polynomial Gauss sums in the special case of $\psi = E$. 

\begin{theorem}[Hayes]\label{thm:hayes}
  If $H$ is a polynomial in $\fnum_q[x]$ with $\deg(H) = m$, then for any multiplicative character $\chi$ of $\fnum_q[x]$, we have 
  \begin{equation}
    \label{eq:hayes}
    (-1)^m G(\chi^{(n)}, E^{(n)}) = \left((-1)^m G(\chi, E) \right)^n.
  \end{equation}
\end{theorem}

The main purpose of this paper is to generalize the above theorem to all of polynomial Gauss sums, we present a compltely analogue of Davenport--Hasse's theorem in polynomial case. 

To state our result, first we note that the character $E_{\lambda}(G, H)$ given by \eqref{eq:Elambda} are all of additive characters $\psi$ modulo $H$ on $\fnum_q[x]$. In other words, for any additive character $\psi$ modulo $H$ on $\fnum_q[x]$, there exists a unique polynomial $G$ in $\fnum_q[x]$, such that $\deg(G)<\deg(H)$, and $\psi = E_{\lambda}(G, H)$ (see Lemma \ref{lem:uniquePoly} below). We write $\psi = \psi_G = E_{\lambda}(G,H)$, and call $G$ the associated polynomial to $\psi$, then we have

\begin{theorem}\label{thm:main}
  If $H$ is a polynomial in $\fnum_q[x]$ of degree $m$, $\chi$ and $\psi$ are multiplicative and additive characters, not both are principal, then we have
  \begin{equation}
    \label{eq:main}
    (-1)^{m-m_1}\frac{\phi^{(n)}(N)}{\phi^{(n)}(H)} G(\chi^{(n)}, \psi^{(n)}) = \left((-1)^{m-m_1} \frac{\phi(N)}{\phi(H)}G(\chi, \psi) \right)^n,
  \end{equation}
where $N=\frac{H}{(G, H)}$, and $G$ is the associated polynomial to $\psi$ ($\psi = \psi_G$), $m_1 = \deg(G, H)$, $\phi(H)$ is the Euler function on $\fnum_q[x]$, and $\phi^{(n)}(H)$ is the function on $\fnum_{q^n}[x]$. In particular, if $H=P^k$, a power of an irreducible $P$, then we have 
\begin{equation}
  \label{eq:mainParticular}
  (-1)^{m-m_1}G(\chi^{(n)}, \psi^{(n)}) = \left( (-1)^{m-m_1}G(\chi, \psi) \right)^n.
\end{equation}
\end{theorem}\vspace{4mm}

If $G=1$ is a constant polynomial, then equality \eqref{eq:main} becomes Hayes' result (Theorem \ref{thm:hayes}). If $(G, H)= 1$, we also have 
\begin{equation}
  \label{eq:mainHayes}
   (-1)^{m}G(\chi^{(n)}, \psi^{(n)}) = \left( (-1)^{m}G(\chi, \psi) \right)^n.
\end{equation}
The equality above essentially belongs to Hayes \cite{H1966a}. 

As we have known, Davenport--Hasse's theorem plays an important role for the rationality of the Zeta function associated to a hypersurface, we wish the result presented here are helpful for the conrugent Zeta function. Finally, we mention a result given by Thakur that an analogue of Davenport--Hasse's theorem for Gauss sums taking values in function fields of one variable over a finite field holds, see Thakur \cite{T1988g}. 

Throughout this paper, by positive polynomial means the polynomial of the leading coefficients is unit in $\fnum_q$, the capital letters $A, B, C,\ldots$ denote polynomials in $\fnum_q[x]$, or $\fnum_{q^n}[x]$, and $a, b, c, \ldots$ denote the elements in $\fnum_q$ or $\fnum_{q^n}[x]$. The absolute value function $|H| = q^m$, where $m=\deg(H)$, and $|H|_n = q^{nm}$ on $\fnum_{q^n}[x]$, which are the numbers of a complete residue class module $H$ on $\fnum_q[x]$ and $\fnum_{q^n}[x]$. 

\section{Properties of character $E_{\lambda}(G, H)$}
\label{sec:prop}
In this section, we first determine the construction of the additive character group modulo $H$ on $\fnum_q[x]$ by using $E_{\lambda}(G, H)$. 

\begin{lemma}\label{lem:uniquePoly}
  For any $\psi$, an additive character modulo $H$ on $\fnum_q[x]$, there exists a unique polynomial $G$ in $\fnum_q[x]$, such that $\deg(G) < \deg(H)$, and $\psi=E_{\lambda}(G,H)$. 
\end{lemma}

\begin{proof}
  For the convenient sake, we write
  \begin{equation}
    \label{eq:psiG}
    \psi_G = E_{\lambda}(G, H). 
  \end{equation}
By \eqref{eq:Elambda}, we have $\psi_{G_1} = \psi_{G_2}$, whenever $G_1 \equiv G_2 (\moda H)$, so we may set $G$ in a complete residue class modulo $H$ in $\fnum_q[x]$, and then $\deg(G)<\deg(H)$. It is easy to show that by \eqref{eq:Elambda}, for any $G_1, G_2$ in $\fnum_q[x]$ that 
\begin{equation}
  \label{eq:psiGbar}
  \psi_{G_1+G_2} = \psi_{G_1}\cdot \psi_{G_2}, \mbox{ and } \bar{\psi}_{G} = \psi_{-G}.
\end{equation}
Since $\psi_G = \psi_0$ the principal character modulo $H$ on $\fnum_q[x]$, if $G=0$, or $H|G$. Conversely, we have $\psi_G=\psi_0$, if and only if $H|G$. Since $\lambda$ is a non-principal character on $\fnum_q$ by assumption, then there is an element $a$ in $\fnum_q$, so that $\lambda(a)\neq 1$. Now if $H\nmid G$, we may let
\begin{equation}
  \label{eq:Rpoly}
  R = (G, H)= a_kx^k + \cdots + a_1x + a_0 \in \fnum_q[x],
\end{equation}
where $0\leq k \leq m -1$, and $a_k\neq 0$. It follow that
$$a\cdot a_k^{-1} x^{m-1-k} R = a x^{m-1}+\cdots .$$
We note that the congruent equation in variable $T$ that 
\begin{equation}
  \label{eq:T}
  GT \equiv a\cdot a_k^{-1} x^{m-1-k} R (\moda H),
\end{equation}
is solvable in $\fnum_q[x]$, therefore, there exists a polynomial $A$ in $\fnum_q[x]$, such that
\begin{equation}
  \label{eq:A}
  GA \equiv a\cdot a_k^{-1} x^{m-1-k} R (\moda H),
\end{equation}
and we see that $t_G(A)=a$ by the definition of \eqref{eq:tG}, and $\psi_G(A) = \lambda(t_G(A)) = \lambda(a) \neq 1$, and $\psi_G \neq \psi_0$. By \eqref{eq:psiGbar}, we have immediately that
\begin{equation}
  \label{eq:psiGneq}
  \psi_{G_1} \neq \psi_{G_2}, \mbox{ if } G_1 \not \equiv G_2 (\moda H). 
\end{equation}
Since if $\psi_{G_1} = \psi_{G_2}$, then $\psi_{G_1-G_2} = \psi_0$, and $G_1\equiv G_2 (\moda H)$. This shows that $\psi_{G}$ are different from each other when $G$ running through a complete residue system of modulo $H$, hence there are exactly $|H| = q^m$ different characters $\psi_G$, but the number of additive characters modulo $H$ on $\fnum_q[x]$ exactly is $q^m$, thus every additive character $\psi$ is just the form of $\psi_G$. We complete the proof of Lemma~\ref{lem:uniquePoly}.
\end{proof}

The next lemma is not new, one may find in Carlitz \cite{C1947t} (see \cite{C1947t}(2.4), (2.5), and (2.6)), but we give a more explicit expression. 

\begin{lemma}\label{lem:carlitz}
  If $A$ is a positive polynomial in $\fnum_q[x]$, then we have
  \begin{equation}
    \label{eq:lem22}
    E_{\lambda}(GA, HA) = E_{\lambda}(G, H).
  \end{equation}
\end{lemma}

\begin{proof}
  For any $B \in \fnum_q[x]$, let 
$$GB \equiv a_{m-1}x^{m-1} + \cdots + a_1x+a_0 (\moda H).$$
Then 
$$AGB \equiv A(a_{m-1}x^{m-1} + \cdots + a_1x+a_0) (\moda AH).$$
Because of $A$ a positive polynomial, we see that the function $t_{GA}$ modulo $HA$ just is the function $t_G$ modulo $H$. Since that $t_{GA}(B) = a_{m-1}$ modulo $HA$, which is just $t_G(B)$ modulo $H$. It follows that
\begin{equation}
  \label{eq:lem22main}
  E_{\lambda}(GH, HA)(B) = \lambda(t_G(B)) = E_{\lambda}(G, H)(B),
\end{equation}
and we have the lemma at once.
\end{proof}

\begin{lemma}\label{lem:twoCases}
  For any $A$ in $\fnum_q[x]$, we have
  \begin{equation*}
    \sum_{G\moda H}\psi_G(A) = \left\{
      \begin{tabular}[l]{l}
        \big|H\big|, \mbox{ if }{H \big| A}\\
        0, \mbox{otherwise}
      \end{tabular}\right.,
  \end{equation*}
where the summation extending over a complete residue system modulo $H$. 
\end{lemma}

\begin{proof}
  By Lemma \ref{lem:uniquePoly}, it just is the orthogonal property of characters. We have the lemma immediately. 
\end{proof}

\section{The separable polynomial Gauss sums}
\label{sec:separable}
The theory for conductors of modulo a polynomial parallels the theory for conductors of Dirichlet characters defined on the integers (see, e.g., \cite[pp. 165--172]{A2013i}), but for proving our theorem, we still state and prove a few basic results. First by Lemma \ref{lem:uniquePoly}, all of the Gauss sums on $\fnum_q[x]$ may be written by
\begin{equation}
  \label{eq:GaussAll}
  G(\chi, \psi) = G(\chi, \psi_G) = \sum_{D\moda H}\chi(D)\psi_G(D),
\end{equation}
where $G\in \fnum_q[x]$ is the polynomial associated to $\psi$. If $(G, H)=1$, it is easy to verify that 
\begin{equation}
  \label{eq:xBar}
  G(\chi, \psi_G) = \bar{\chi}(G) G(\chi, \psi_1).
\end{equation}

\begin{definition}\label{def:sep}
  A Gauss sum $G(\chi, \psi_G)$ is said to be separable if $G(\chi, \psi_G) = \bar{\chi}(G) G(\chi, \psi_1)$. 
\end{definition}

By \eqref{eq:xBar}, if $(G, H)=1$, then $G(\chi, \psi_G)$ is separable. For the case of $(G, H)>1$, then $G(\chi, \psi_G)$ is separable if and only if $G(\chi, \psi_G) = 0$. The following lemma gives an important consequence of separability. 

\begin{lemma}\label{lem:sep}
  If $G(\chi, \psi_G)$ is separable for every $G$ in $\fnum_q[x]$, then
  \begin{equation}
    \label{eq:lemSep}
    |G(\chi, \psi_1)|^2 = |H| = q^m.
  \end{equation}
\end{lemma}

\begin{proof}
  \begin{equation}
    \label{eq:lemSepProof}
    \begin{split}
|G(\chi, \psi_1)|^2 &= G(\chi, \psi_1)\overline{G(\chi, \psi_1)}\\
&= G(\chi, \psi_1)\sum_{D\moda H} \bar{\chi}(D)\psi_1(-D)\\
&= \sum_{D\moda H}G(\chi, \psi_D)\psi_D(-1)\\
&= \sum_{A\moda H}\chi(A)\sum_{D\moda H}\psi_D(A-1).
    \end{split}
  \end{equation}
The inner sum in above is zero by Lemma \ref{lem:twoCases}, if $A\not \equiv 1 (\moda H)$, so we have Lemma \ref{lem:sep}.
\end{proof}

\begin{lemma}\label{lem:existence}
  If $G(\chi, \psi_G)\neq 0$ for some $G$ in $\fnum_q[x]$ with $(G, H)>1$, then there exists a polynomial $N$ in $\fnum_q[x]$, such that $N|H$, $\deg(N)<\deg(H)$, and
  \begin{equation}
    \label{eq:lemExist}
    \chi(A)=1, \mbox{ whenever }(A, H)=1, \mbox{ and }A\equiv 1 (\moda N). 
  \end{equation}
\end{lemma}

\begin{proof}
  For given $G$, and $G(\chi, \psi_G)\neq 0$, $(G, H)>1$. Let $N=H\cdot (G, H)^{-1}$, thus $N|H$, and $\deg(N)<\deg(H)$. If $(A, H)=1$, then
  \begin{equation}\label{eq:lemExistProof}
    \begin{split}
      G(\chi, \psi_G) &= \sum_{D \moda H}\chi(AD)\psi_G(AD)\\
&= \chi(A)\sum_{D \moda H}\chi(D)\psi_G(AD).
    \end{split}
  \end{equation}
If $A \equiv 1 (\moda N)$, we write
$$A = 1+BN, \mbox{ for some }B\in \fnum_q[x],$$
and then
$$AG = G+BNG = G+BHG(H, G)^{-1}.$$
So we have $AG \equiv G (\moda H)$, and $\psi_G(AD) = \psi_{GA}(D) = \psi_G(D)$. Therefore, equation \eqref{eq:lemExistProof} becomes that
\begin{equation}
  \label{eq:lemExistProof2}
  G(\chi, \psi_G) = \chi(A) \sum_{D \moda H}\chi(D)\psi_G(D),
\end{equation}
and we have $\chi(A)=1$ because of $G(\chi, \psi_G)\neq 0$. We complete the proof of Lemma \ref{lem:existence}.
\end{proof}

\begin{definition}
  A polynomial $N$ in $\fnum_q[x]$ is called an induced modulu of $\chi$ if $N|H$, and for $(A,H)=1$
  \begin{equation}
    \label{eq:defInduce}
    \chi(A)=1, \mbox{ whenever } A\equiv 1 (\moda N).
  \end{equation}
\end{definition}

By the definition, we see that $H$ itself is an induced modulu of any $\chi$. Moreover, as a direct consequence of Lemma \ref{lem:existence}, we also have

\begin{corollary}\label{cor:existence}
  If $(G, H)>1$, and $G(\chi, \psi_G)\neq 0$, then $N=H(G,H)^{-1}$ is an induced modulu of $\chi$. 
\end{corollary}

\begin{lemma}\label{lem:Equiv}
Let $N|H$, then $N$ is an induced modulu of $\chi$ if and only if for any $A,B$ in $\fnum_q[x]$, and $(AB, H)=1$, we have
  \begin{equation}
    \label{eq:lemEquiv}
    \chi(A)=\chi(B), \mbox{ whenever }A\equiv B (\moda N).
  \end{equation}
\end{lemma}

\begin{proof}
  If \eqref{eq:lemEquiv} holds, let $B=1$, then $N$ is an induced modulu of $\chi$. Conversely, if $N$ is an induced modulu of $\chi$, suppose $(A, H)=(B, H)=1$, $A\equiv B (\moda N)$, and let $B\cdot B^{-1}\equiv 1 (\moda H)$. Then
$BB^{-1}\equiv 1 (\moda N)$, and $AB^{-1}\equiv 1 (\moda N)$. Thus $\chi(AB^{-1})=1$, and $\chi(A)=\chi(B)$, the lemma follows. 
\end{proof}

\begin{lemma}\label{lem:multi}
  If $N|H$, and $N$ is an induced modulu of $\chi$, then $\chi$ can be expressed as a product
  \begin{equation}
    \label{eq:lemMulti}
    \chi = \chi_0 \delta,
  \end{equation}
where $\chi_0$ is the principal multiplicative modulo $H$, and $\delta$ is a multiplicative character modulo $N$. 
\end{lemma}

\begin{proof}
  If $\chi = \chi_0 \delta$, trivially, $N$ is an induced modulu of $\chi$. Conversely, if $N$ is an induced modulu, we may determine a character $\delta$ modulo $N$ by setting $\delta(A)=0$, if $(A,N)>1$. If $(A,N)=1$, one may find a polynomial $B$ in $\fnum_q[x]$, so that
  \begin{equation}
    \label{eq:lemMultiProof}
    (B, H)=1, \mbox{ and }B \equiv A (\moda N).
  \end{equation}

Since the arithmetic progress $\{A+RN | R\in \fnum_q[x]\}$ contains infinitely many irreducibles (see Artin \cite{A1924q}, or \cite[Theorem~4.7]{R2002n}), so we may choose one that does not divide $H$ and call this $B$, which is unique modulo $N$ clearly. Now we define $\delta(A) = \chi(B)$. The number $\delta(A)$ is well-defined because $\chi$ takes equal values at polynomials which are congruent modulo $N$ and relatively prime to $H$. By this determination, we can easily verify that $\delta$ is, indeed, a character modulo $N$, and \eqref{eq:lemMultiProof} holds. This is the proof of Lemma \ref{lem:multi}.
\end{proof}

\begin{definition}
  An induced modulu $N$ of $\chi$ is called the conductor of $\chi$, if $N$ is positive, and an induced modulu of $\chi$, and the degree of $N$ is least among all of induced modulus of $\chi$. We denote by $C_{\chi}$ the conductor of $\chi$. If $C_{\chi} = H$, then we call $\chi$ is a primitive character. 
\end{definition}

As a consequence of Lemma \ref{lem:multi}, we have

\begin{corollary}
  If $C_{\chi}$ is the conductor of $\chi$, then $\chi$ can be expressed as a product $\chi = \chi_0\delta$, where $\delta$ is a primitive character modulo $C_{\chi}$. 
\end{corollary}

The following lemma is well-known that
\begin{lemma}\label{lem:hayesCon}
  The conductor of $\chi$ divides every induced modulu of $\chi$. 
\end{lemma}

\begin{proof}
  See Hayes \cite[Theorem~4.2]{H1966a}.
\end{proof}

We have an alternate description of primitive character as the case of Dirichlet characters.

\begin{lemma}\label{lem:primitveSep}
  Let $\chi$ be a character modulo $H$, then $\chi$ is primitive, if and only if the Gauss sums $G(\chi, \psi_G)$ is separable for every polynomial $G$. 
\end{lemma}

\begin{proof}
   If $\chi$ is primitive, then $G(\chi, \psi_G)$ is separable by Lemma \ref{lem:existence}, so we only prove the converse. It suffices to prove that if $\chi$ is not primitive, then there exists some $G$ with $(G,H)>1$, and $G(\chi, \psi_G)\neq 0$. If $\chi$ is not primitive, let $C_{\chi}$ be the conductor of $\chi$, $N=\frac{H}{C_{\chi}}$, then $(N,H)>1$ by $\deg(C_{\chi}) < \deg(H)$, moreover $G(\chi, \psi_N)\neq 0$. Since $\chi = \chi_0\delta$, where $\delta$ is primitive character modulo $C_{\chi}$. By Lemma \ref{lem:carlitz}, we have
    \begin{equation}
      \label{eq:lemPrimitiveSep}
      \begin{split}
        G(\chi, \psi_N) &= \sum_{D \moda H} \chi_0(D)\delta(D)\psi_N(D)\\
&= \sum_{\substack{D \moda H, \\(D,H)=1}}\delta(D)E_{\lambda}(N, H)(D)\\
&= \sum_{\substack{D \moda H,\\ (D,H)=1}}\delta(D)E_{\lambda}(1, C_{\chi})(D)\\
&= \frac{\phi(H)}{\phi(C_{\chi})} \sum_{D\moda C_{\chi}}\delta(D)\psi_1(D)\\
&= \frac{\phi(H)}{\phi(C_{\chi})} G(\delta, \psi_1),
      \end{split}
    \end{equation}
where $G(\delta, \psi_1)$ is a Gauss sum modulo $C_{\chi}$. By Lemma \ref{lem:sep}
\begin{equation}
  \label{eq:lemPrimitiveSep2}
  |G(\delta, \psi_1)|^2 = |C_{\chi}|. 
\end{equation}
So we have $G(\chi, \psi_N)\neq 0$, and the lemma follows. 
\end{proof}

\section{Separable Gauss sums on $\fnum_{q^n}[x]$}
\label{sec:sepExt}

For a polynomial $H$ in $\fnum_q[x]$, and also a polynomial in $\fnum_{q^n}[x]$. Let $\psi_G^{(n)}$ be an additive character, $\chi^{(n)}$ a multiplicative character modulo $H$ on $\fnum_{q^n}[x]$ given by \eqref{eq:gaussSum}. We recall that $\psi_G = E_{\lambda}(G, H)$ in $\fnum_q[x]$, it still holds in $\fnum_{q^n}[x]$, namely
\begin{equation}
  \label{eq:ElambdaFqn}
  \psi_G^{(n)} = E_{\lambda}^{(n)}(G, H),\mbox{ if }G\in \fnum_q[x].
\end{equation}
Since for any polynomial $A$ in $\fnum_{q^n}[x]$
\begin{equation*}
  E_{\lambda}^{(n)}(G, H)(A) = E_{\lambda}(G, H)(\tr(A))=\lambda(t_G(\tr(A)))=\psi_G^{(n)}(A).
\end{equation*}
The last equality is because of $\tr(G\cdot A) = G\tr(A)$ for $G\in \fnum_q[x]$. So we have if $G\in \fnum_q[x]$ then
\begin{equation}
  \label{eq:psiN}
  \psi_G^{(n)}(A) = \psi_1^{n}(GA).
\end{equation}

Now, Lemma \ref{lem:carlitz} becomes that

\begin{lemma}\label{lem:carlitzRe}
  If $A\in \fnum_q[x]$ is a positive polynomial, and $G\in \fnum_q[x]$, then
  \begin{equation}
    \label{eq:carlitzRe}
    E_{\lambda}^{(n)}(GA, HA) = E_{\lambda}^{(n)}(G, H).
  \end{equation}
\end{lemma}

\begin{proof}
  For any $B\in \fnum_{q^n}[x]$, let
$$G\tr(B) \equiv a_{m-1}x^{m-1} + \cdots + a_1x+a_0 (\moda H).$$
Then 
$$AG\tr(B) \equiv A(a_{m-1}x^{m-1} + \cdots + a_1x+a_0) (\moda AH).$$
By definition of $E_{\lambda}(G, H)$ and $E_{\lambda}^{(n)}(G, H)$, it follows that
$$E_{\lambda}^{(n)}(GA, HA)(B) = E_{\lambda}(GA, HA)(\tr(B)) = \lambda(t_{GA}(\tr(B))) = \lambda(a_{m-1}),$$
and
$$E_{\lambda}^{(n)}(G, H)(B) = E_{\lambda}(G, H)(\tr(B)) = \lambda(t_{G}(\tr(B))) = \lambda(a_{m-1}).$$
This lemma follows at once.
\end{proof}

\begin{lemma}\label{lem:inducedFqn}
  If $N\in \fnum_q[x]$ is an induced modulu of $\chi$, then $N$ is also an induced modulu of $\chi^{(n)}$ on $\fnum_{q^n}[x]$. 
\end{lemma}

\begin{proof}
  Suppose $A\in \fnum_{q^n}[x]$, $(A, H)=1$, and $A = 1 (\moda N)$, it is easily seen that if $N(A)$ is the norm that
  \begin{equation}
    \label{eq:lemInducedFqn}
    N(A)\equiv 1 (\moda N).
  \end{equation}
This equality is the Theorem~2.1 of Hayes \cite{H1966a}. Then $\chi^{(n)}(A)=\chi(N(A))=1$, and $N$ is an induced modulu of $\chi^{(n)}$ on $\fnum_q[x]$. 
\end{proof}

The following lemmas is due to Hayes \cite{H1966a} that
\begin{lemma}\label{lem:HayesConductor}
  Let $C_{\chi}$ be the conductor of $\chi$ on $\fnum_q[x]$, and $C_{\chi^{(n)}}$ be the conductor of $\chi^{(n)}$ on $\fnum_{q^n}[x]$, then we have $C_{\chi} = C_{\chi^{(n)}}$.
\end{lemma}

\begin{proof}
  See Hayes \cite[Theorem~4.5]{H1966a}.
\end{proof}

As a direct corollary of the above lemma and Lemma \ref{lem:primitveSep}, we have
\begin{corollary}
  Let $G(\chi, \psi)$ be the Gauss sums modulo $H$ on $\fnum_q[x]$, $G(\chi^{(n)}, \psi^{(n)})$ the Gauss sums modulo $H$ on $\fnum_{q^n}[x]$, then $G(\chi^{(n)}, \psi^{(n)})$ is separable if and only if $G(\chi, \psi)$ is separable. 
\end{corollary}

\begin{lemma}\label{lem:inducedBoth}
  Suppose $G$ and $H$ are in $\fnum_q[x]$, and $(G, H)>1$, if $G(\chi, \psi_G) \neq 0$, or $G(\chi^{(n)}, \psi_G^{(n)})\neq 0$, then $N = \frac{H}{(G, H)}$ is an induced modulu of both $\chi$ and $\chi^{(n)}$. 
\end{lemma}

\begin{proof}
  It is easily seen that if $N$ is an induced modulu of $\chi$, then any multiple of $N$ which divides $H$, again is an induced modulu of $\chi$. To prove the lemma, first suppose $G(\chi, \psi_G)\neq 0$, because $(G,H)>1$, then by Lemma \ref{lem:existence} and Corollary \ref{cor:existence}, $N$ is an induced modulu of $\chi$. By Lemma \ref{lem:inducedFqn}, $N$ is also an induced modulu of $\chi^{(n)}$, the lemma holds. If $G(\chi^{(n)},\psi_G^{(n)})\neq 0$, then $N$ is an induced modulu of $\chi^{(n)}$. Let $C_{\chi^{(n)}}$ be the conductor of $\chi^{(n)}$, then by Lemma \ref{lem:hayesCon}, we have $C_{\chi^{(n)}} | N$, and then $C_{\chi}|N$, where $C_{\chi}$ is the conductor of $\chi$. Therefore $N$ is an induced modulu of $\chi$. We complete the proof of Lemma \ref{lem:inducedBoth}.
\end{proof}

\section{Proof of Theorem \ref{thm:main}}
\label{sec:mainProof}
We consider two cases to prove this theorem. First, if $(G,H)=1$, then $(G, H)=1$ in $\fnum_{q^n}[x]$. The Gauss sums on $\fnum_{q^n}[x]$ is that
\begin{equation}
  \label{eq:proofMain1}
  G(\chi^{(n)},\psi_G^{(n)}) = \bar{\chi}^{(n)}(G)G(\chi^{(n)},\psi_1^{(n)}),
\end{equation}
and the Gauss sums $G(\chi, \psi_G)$ on $\fnum_q[x]$ is following
\begin{equation}
  \label{eq:proofMain2}
  G(\chi, \psi_G) = \bar{\chi}(G)G(\chi, \psi_1).
\end{equation}
We note that $\chi^{(n)}(G) = \chi^n(G)$, $\psi_1=E$, by Theorem \ref{thm:hayes}, we have
\begin{equation}
  \label{eq:proofMain3}
  \begin{split}
 (-1)^mG(\chi^{(n)},\psi_G^{(n)}) &= (-1)^m\bar{\chi}^n(G)G(\chi^{(n)},\psi_1^{(n)})\\
&= \left( (-1)^m \bar{\chi}(G)G(\chi, \psi_1) \right)^n \\
&= \left( (-1)^m G(\chi, \psi_G) \right)^n.
  \end{split}
\end{equation}
Because of $(G, H)=1$, then $m_1=\deg(G,H)=0$, and $\phi^{(n)}(H)=\phi^{(n)}(N)$, $\phi(H)=\phi(N)$, \eqref{eq:main} of Theorem \ref{thm:main} holds in the case of $(G,H)=1$. 

Next, we suppose $(G,H)>1$, and let $H_1 = \frac{H}{(G, H)}$, $G_1 = \frac{G}{(G,H)}$, thus $(G_1, H_1) = 1$. If both of $G(\chi, \psi_G)$ and $G(\chi^{(n)}, \psi_G^{(n)})$ are zero, then \eqref{eq:main} is trivial. Therefore, we may assume $G(\chi, \psi_G) \neq 0$, or $G(\chi^{(n)},\psi_G^{(n)})\neq 0$. By this assumption, then $H_1$ is an induced modulu of both $\chi$ and $\chi^{(n)}$. By Lemma \ref{lem:Equiv}, we may write
\begin{equation}
  \label{eq:proofMain4}
  \chi = \chi_0\delta, \mbox{ and }\chi^{(n)}=\chi_0^{(n)}\delta^{(n)},
\end{equation}
where $\delta$ is a multiplicative character modulo $H_1$, and $\delta^{(n)}(A) = \delta(N(A))$  is a multiplicative character modulo $H_1$ on $\fnum_{q^n}[x]$. The Gauss sums on $\fnum_{q^n}[x]$ is that
\begin{equation}
  \label{eq:proofMain5}
  \begin{split}
  G(\chi^{(n)},\psi_G^{(n)})&=\sum_{\substack{D\moda H,\\ D\in \fnum_{q^n}[x]}} \chi_0^{(n)}(D)\delta^{(n)}(D)\psi_G^{(n)}(D)\\
&= {\sum_{\substack{D\moda H,\\ D\in \fnum_{q^n}[x]}}}' \delta^{(n)}(D)E_{\lambda}^{(n)}(G, H)(D)\\
&= {\sum_{\substack{D\moda H,\\ D\in \fnum_{q^n}[x]}}}' \delta^{(n)}(D)E_{\lambda}^{(n)}(G_1, H_1)(D)\\
&= \frac{\phi^{(n)}(H)}{\phi^{(n)}(H_1)}G(\delta^{(n)}, \psi_{G_1}^{(n)}),
  \end{split}
\end{equation}
where summation ${\sum}'$ means $(D, H)=1$, and $G(\delta^{(n)},\psi_{G_1}^{(n)})$ is a Gauss sums modulo $H_1$. The Gauss sums $G(\chi, \psi_G)$ modulo $H$ on $\fnum_q[x]$ is that
\begin{equation}
  \label{eq:proofMain6}
  \begin{split}
    G(\chi, \psi_G) &= {\sum_{D\moda H}}' \delta(D)\psi_G(D)\\
&= {\sum_{D\moda H}}' \delta(D)E_{\lambda}(G,H)(D)\\
&= {\sum_{D\moda H}}' \delta(D)E_{\lambda}(G_1,H_1)\\
&= \frac{\phi(H)}{\phi(H_1)}G(\delta, \psi_{G_1}),
  \end{split}
\end{equation}
where $G(\delta, \psi_{G_1})$ is a Gauss sums modulo $H_1$. Because of $(G_1, H_1)=1$, the discussion for first case gives us that
\begin{equation}
  \label{eq:proofMain7}
  (-1)^{m-m_1}G(\delta^{(n)},\psi_{G_1}^{(n)}) = \left( (-1)^{m-m_1}G(\delta, \psi_{G_1}) \right)^n,
\end{equation}
where $m-m_1 = \deg(H_1)$, and the equality \eqref{eq:main} of Theorem \ref{thm:main} follows immediately. 

To prove \eqref{eq:mainParticular} of Theorem \ref{thm:main}, if $H|G$, then $\psi_G = \psi_0$ is the principal character modulo $H$, then both sides of \eqref{eq:mainParticular} are zero, if $\chi$ is non-principal, then we may suppose that $H \nmid G$. 

Since $H=P^k$, where $k\geq 1$, and $P$ is an irreducible in $\fnum_q[x]$, it is well-known that $P$ is product of exactly $(h, n)$ irreducibles in $\fnum_{q^n}[x]$, where $h=\deg(P)$ (see \cite[Theorem~2.1]{H1965t}, for example). If $H\nmid G$, then $\frac{H}{(G, H)} = N = P^{k_1}$, where $1\leq k_1 \leq k$, it is easy to verify that
\begin{equation}
  \label{eq:proofMain8}
  \phi^{(n)}(N)(\phi^{(n)}(H))^{-1} = (\phi(N)\phi^{-1}(H))^n.
\end{equation}

So \eqref{eq:mainParticular} follows from \eqref{eq:main}, we complete the proof of Theorem \ref{thm:main}. 

{\small
\bibliographystyle{siam}
\bibliography{GaussSum}

\begin{thebibliography}{10}

\bibitem{A2013i}
{\sc T.~Apostol}, {\em Introduction to Analytic Number Theory},
  Springer--Verlag, New York Heidelberg Berlin, 1976.

\bibitem{A1924q}
{\sc E.~Artin}, {\em {Quadratische K{\"o}rper im Gebiete der h{\"o}heren
  Kongruenzen. II}}, Math. Z., 19 (1924), pp.~207--246.

\bibitem{C1947t}
{\sc L.~Carlitz}, {\em The singular series for sums of squares of polynomials},
  Duke Math. J., 14 (1947), pp.~1105--1120.

\bibitem{DH1934D}
{\sc H.~Davenport and H.~Hasse}, {\em {Die Nullstellen der
  Kongruenzzetafunktionen in Gewissen Zyklischen F{\"a}llen}}, J. Reine. Angew.
  Math., 172 (1934), pp.~151--182.

\bibitem{H1965t}
{\sc D.~Hayes}, {\em The distribution of irreducibles in {$GF[q, x]$}}, Trans.
  Amer. Math. Soc., 117 (1965), pp.~1017--1033.

\bibitem{H1966a}
\leavevmode\vrule height 2pt depth -1.6pt width 23pt, {\em A polynomial
  generalized {Gauss} sum.}, J. Reine. Angew. Math., 222 (1966), pp.~113--119.

\bibitem{IR1982a}
{\sc K.~Ireland and M.~Rosen}, {\em {A Classical Introduction to Modern Number
  Theory}}, vol.~84 of GTM, Springer--Verlag, New York Heidelberg Berlin, 1982.

\bibitem{R2002n}
{\sc M.~Rosen}, {\em {Number Theory in Function Fields}}, vol.~210 of GTM,
  Springer--Verlag, New York Heidelberg Berlin, 2002.

\bibitem{S1937r}
{\sc H.~Schmid}, {\em {Relationen zwischen verallgemeinerten Gau{\ss}schen
  Summen.}}, J. Reine. Angew. Math., 176 (1937), pp.~189--191.

\bibitem{T1988g}
{\sc D.~Thakur}, {\em Gauss sums for {$F_q[T]$}}, Invent. Math., 94 (1988),
  pp.~105--112.

\end{thebibliography}
}
\end{document}